\documentclass[12pt]{article}
\usepackage[T1]{fontenc}
\usepackage[utf8]{inputenc}
\usepackage[a4paper, margin=2.7cm]{geometry}
\usepackage{amssymb}
\usepackage{amsmath}
\usepackage{graphicx}
\usepackage{hyperref}
\usepackage{url}
\newcommand{\diff}[0]{\mathrm{d} }
\newcommand{\rea}[0]{\mathrm{Re} }
\newcommand{\ima}[0]{\mathrm{Im} }

\newcommand{\teg}[2]{\int\limits_{#1}^{#2}   }
\newcommand{\cont}[1]{\oint\limits_{#1} }
\newcommand{\nsum}[2]{\sum_{#1}^{#2}}
\newcommand{\resi}[1]{\underset{#1}{\textup{Res}}}

\begin{document}
\begin{center}
      {\bf  New proof of the Gaussian integral using the residue\\ theorem with links to the Riemann Zeta function}\\[5mm]
  Bastien Jean Quemener
\end{center}
\vskip+0.6cm
\begin{abstract}
    In this paper the Gaussian integral is proven using contour integration on $\frac{1}{e^{z^2}+1}$ and linking it using a limit to said Gaussian integral. The limit is also related to the Riemann Zeta function using a few manipulations. This new and original proof comes as an addition to the already many pre-existing proofs of the Gaussian integral \cite{conrad}
\end{abstract}
\vskip+0.2cm \noindent  
\tableofcontents
\section{The proof}
Let us find $\teg{-\infty}{\infty} \frac{1}{e^{x^2}+1}\diff x$, because $z^2$ and $e^z$ are holomorphic, and $1/(z+1)$ is also holomorphic, we get that $1/(e^{z^2}+1)$ is holomorphic, hence we can use contour integration for our integral.
\newline
The poles of $1/(e^{z^2}+1)$ at $z=c$ correspond to the values of $c$ such that $e^{c^2}+1=0$, $e^{c^2}=-1$ :
\vskip-0.5cm
$$e^{c^2}=e^{(2k+1)\pi i}\;\;\;\; \;\;\;(k\in \mathbb{N})$$
$${c^2}={(2k+1)\pi i}=(2k+1)\pi e^{i \pi/2}$$
With $k$ running from $-\infty$ to $\infty$, $2k+1$ can be reduced to $\pm(2k+1)$ with $k$ running from $0$ to $\infty$
$${c^2}=\pm(2k+1)\pi e^{i \pi/2}\;\;\;\; \;\;\;(k\in \mathbb{N}, 0\le k)$$
$$c=\pm\sqrt{\pm}\sqrt{(2k+1)\pi} e^{i \pi/4}=\pm\sqrt{\pm}\sqrt{(2k+1)\pi} \frac{1+i}{\sqrt{2}}=\sqrt{(2k+1)\pi} \frac{\pm 1+\pm i}{\sqrt{2}}$$
Our contour is going to be in the upper complex plane, hence poles with negative imaginary parts are irrelevant, all our relevant poles are then :
\begin{equation}\label{poles}
  \boxed{A_k=( i+1)\sqrt{\frac{(2k+1)\pi}{2}}\,\,\,( 0\le k, k\in \mathbb{N})}
\end{equation}
$$  \boxed{B_k=( i-1)\sqrt{\frac{(2k+1)\pi}{2}}\,\,\,( 0\le k, k\in \mathbb{N})}$$
\newline
$H_n$ is a real value dependant of $n$. We take our contour, $G_n$, to be a line starting from $-H_n$ to $H_n$, then, a semi circle with positive imaginary value centered at 0 of radius $H_n$ starting at the point $H_n$ and ending back up at $-H_n$, called $F_n$. By taking $H_n$ to be between the absolute value of $A_n$ and $A_{n-1}$ or $B_n$ and $B_{n-1}$, this ensures we are not treading on any poles and that all the poles in our contour are $A_k$ and $B_k$ with $k$ going from $0$ to $n-1$.
\begin{figure}[!ht]
  \centering
  \caption{The contour $G$ in the complex plane}
  \includegraphics[scale=0.3]{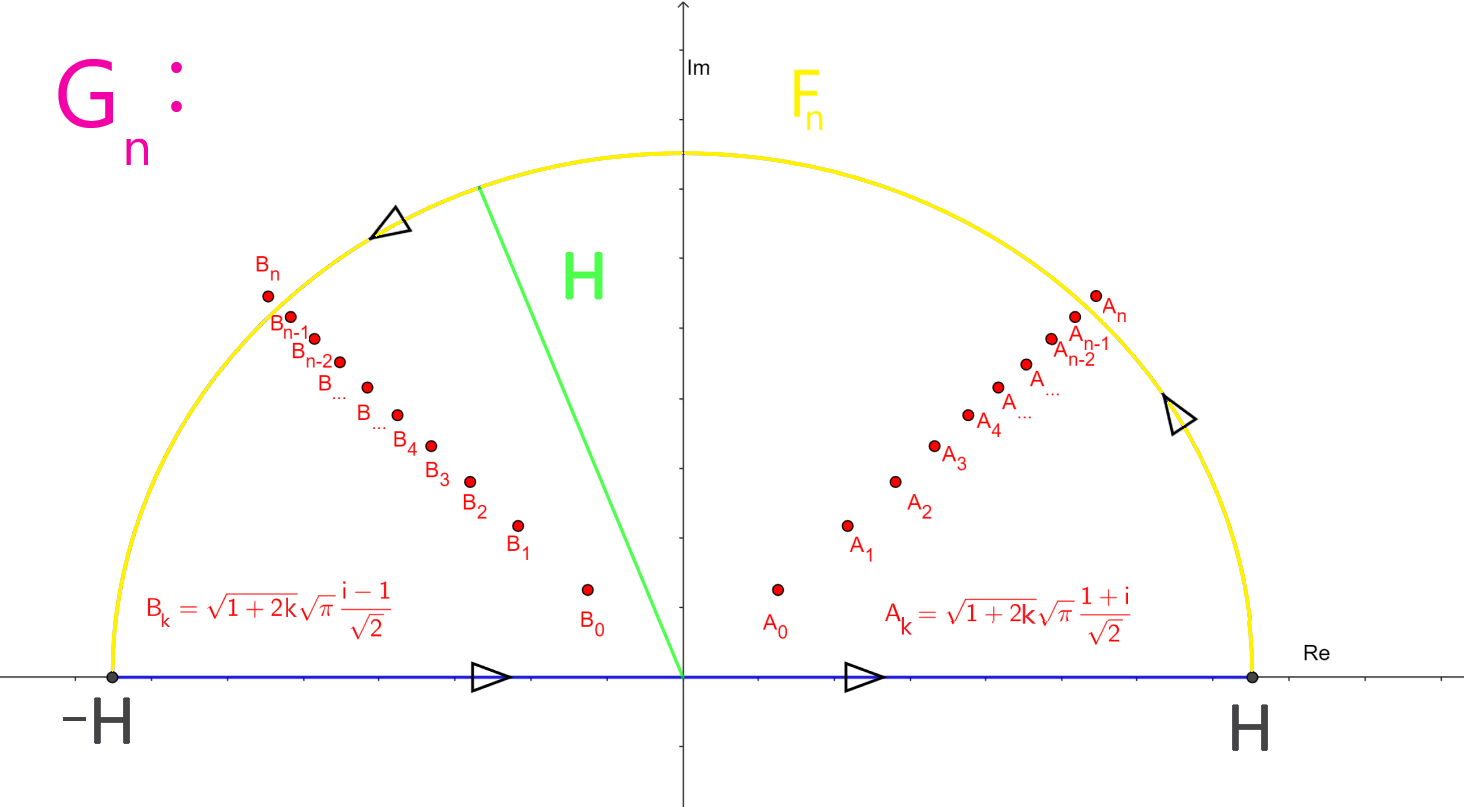}
\end{figure}
Integrating over $G_n$ gives:
$$\cont{G_n} \frac{1}{e^{z^2}+1} \diff z = \teg{-H_n}{H_n}\frac{1}{e^{x^2}+1} \diff x+\teg{F_n}{}\frac{1}{e^{z^2}+1} \diff z$$
From our contour, $F_n$ can be parametrized as :
$$F_n : z=H_n e^{iv} \;\;\;\;\;\;\;\;\;\;\;\;\;\;\;\;\;\;\;\;\; 0 \to v \to \pi   \;\;\;\;\;\;\;\;\;\;\;\;\;\;\;\;\;\;\;\;\;\diff z = H_n i e^{iv} \diff v$$
\begin{equation} \label{norescont}
  \boxed{\cont{G_n} \frac{1}{e^{z^2}+1} \diff z = \teg{-H_n}{H_n}\frac{1}{e^{x^2}+1} \diff x+\teg{0}{\pi }\frac{H_n i e^{iv} }{e^{H_n ^2 e^{2iv}}+1} \diff v}
\end{equation}
Because the poles of $1/(e^{z^2}+1)$ are of first order, the residue becomes :
$$\resi{z=c} \frac{1}{e^{z^2}+1}=\lim_{z \to c} \frac{z-c}{e^{z^2}+1}$$
Since both sides of the fraction go to zero, $\lim_{z\to c} \frac{1}{2z e^{z^2}} $ exists ($c$ is nonzero) and because $\frac{1}{e^{z^2}+1}$ is holomorphic, we can use L'Hopital's rule, such as :
$$\resi{z =c} \frac{1}{e^{z^2}+1} =\lim_{z\to c} \frac{\frac{\diff(z-c)}{\diff z} }{\frac{\diff( e^{z^2}+1)}{\diff z}}= \lim_{z\to c} \frac{1}{2z e^{z^2}}=  \frac{1}{2c e^{c^2}}=\frac{-1}{2c}$$
$H_n$ is between the absolute value of $A_n$ and $A_{n-1}$
$$|A_{n-1}|=\sqrt{(2n-1)\pi}<H_n<\sqrt{(2n+1)\pi}=|A_{n}|$$
At some point, $k=n-1$, because $0\le k, k\in \mathbb{N}$ this makes $1 \le n, n \in \mathbb{N}$. From the residue theorem, because we have only $A_k$ and $B_k$ poles with $k$ from $0$ to $n-1$ in our contour, from  \ref{norescont},  \ref{poles} :
$$\cont{G_n} \frac{1}{e^{z^2}+1}\diff z =2\pi i \nsum{k=0}{n-1} (\resi{z=A_k}\frac{1}{e^{z^2}+1} + \resi{z=B_k}\frac{1}{e^{z^2}+1}) =2\pi i \nsum{k=0}{n-1} (\frac{-1}{2 A_k}-\frac{1}{2B_k}) $$
$$ =2\pi i \nsum{k=0}{n-1} (\frac{-1}{2 ( i+1)\sqrt{\frac{(2k+1)\pi}{2}}}-\frac{1}{2( i-1)\sqrt{\frac{(2k+1)\pi}{2}}}) =\nsum{k=0}{n-1}\frac{-\sqrt{2\pi}}{\sqrt{2k+1}}$$
\vskip-0.3cm
\begin{equation} \label{final no hn}
  \boxed{\teg{-H_n}{H_n}\frac{\diff x}{e^{x^2}+1} +\teg{0}{\pi }\frac{H_n i e^{iv} \diff v}{e^{H_n ^2 e^{2iv}}+1} =\nsum{k=0}{n-1}\frac{-\sqrt{2\pi} }{\sqrt{2k+1}}} 
\end{equation}
$$\boxed{\sqrt{(2n-1)\pi}<H_n<\sqrt{(2n+1)\pi}, \,\,\,\,1\le n, n\in \mathbb{N}}$$
Setting $\sqrt{(2n-1)\pi}< H_n=\sqrt{2n\pi}<\sqrt{(2n+1)\pi}  $ :
\begin{equation}\label{final yes hn complex}
   \boxed{\teg{-\sqrt{2n\pi}}{\sqrt{2n\pi}}\frac{1}{e^{x^2}+1} \diff x+\nsum{k=0}{n-1}\frac{\sqrt{2\pi}}{\sqrt{2k+1}}=-\teg{0}{\pi }\frac{i e^{iv} \sqrt{2n\pi}  }{e^{2n\pi  e^{2iv}}+1} \diff v\;\;\;\;\;\;\;\;\;\;\;\;\;\; ( 1\le n, n\in \mathbb{N})}
\end{equation}
Applying the real part to the equation does not change it because $\teg{-\sqrt{2n\pi}}{\sqrt{2n\pi}}\frac{ \diff x}{e^{x^2}+1}$ and $\nsum{k=0}{n-1}\frac{\sqrt{2\pi}}{\sqrt{2k+1}}$ are real the equation is all real. Because $\rea$ commutes with integration, and that $\rea(iz)=-\ima(z)$:
\vskip-0.3cm
$$-\teg{0}{\pi }\frac{i e^{iv} \sqrt{2n\pi}  }{e^{2n\pi  e^{2iv}}+1} \diff v=-\rea(\teg{0}{\pi }\frac{i e^{iv} \sqrt{2n\pi}  }{e^{2n\pi  e^{2iv}}+1} \diff v)=-\teg{0}{\pi }\rea(\frac{i e^{iv} \sqrt{2n\pi}  }{e^{2n\pi  e^{2iv}}+1} )\diff v$$
We then have, with defining $P(n)$
\begin{equation}\label{final yes hn real}
   \boxed{P(n)=\!\!\!\!\!\!\teg{-\sqrt{2n\pi}}{\sqrt{2n\pi}}\frac{\diff x}{e^{x^2}+1} +\nsum{k=0}{n-1}\frac{\sqrt{2\pi}}{\sqrt{2k+1}}=\sqrt{n}\teg{0}{\pi }\ima(\frac{ e^{iv} \sqrt{2\pi}  }{e^{2n\pi  e^{2iv}}+1} )\diff v\;\;\;\;\;\;\;\; ( 1\le n, n\in \mathbb{N})}
\end{equation}
Because the function $\frac{1}{e^{x^2}+1}$ is positive, the greater our interval of integration is, the greater the integral of the function is, as $n$ strictly increases, the interval from $-\sqrt{2n\pi}$ to $\sqrt{2n\pi}$ strictly increases, hence $P(n)= \teg{-\sqrt{2\pi}}{\sqrt{2n\pi}}\frac{\diff x}{e^{x^2}+1} $ is a strictly increasing and positive function.
\newline
From $1+x\le e^x$ we can have :
\vskip-1cm
$$\frac{1}{e^{x^2}+1} \le \frac{1}{x^2+2}$$
\begin{equation}\label{bounded}
\boxed{0\le P(n)=\teg{-\sqrt{2n\pi}}{\sqrt{2n\pi}} \frac{\diff x}{e^{x^2}+2}\le \teg{-\sqrt{2n\pi}}{\sqrt{2n\pi}} \frac{\diff x}{x^2+2}= \teg{-\sqrt{n\pi}}{\sqrt{n\pi}} \frac{1/\sqrt{2}}{x^2+1}\diff x=\frac{2\arctan(\sqrt{n\pi})}{\sqrt{2}}\le \frac{\pi}{\sqrt{2}}}
\end{equation}
This makes $P(n)$ bounded between $0$ and $\pi/\sqrt{2}$, from the monotone convergence theorem \cite{moncontheo}, if a sequence is increasing and bounded, it's limit is the supremum, applying it to $P(n)$ \ref{final yes hn real}, remembering that $n$ is strictly an integer :
$$ \lim_{n\to\infty} \teg{-\sqrt{2n\pi}}{\sqrt{2n\pi}}\frac{\diff x}{e^{x^2}+1} = \teg{-\infty}{\infty}\frac{ \diff x}{e^{x^2}+1}=\limsup_{n\to\infty}(\sqrt{n} \teg{0}{\pi } \ima(\frac{\sqrt{2\pi}e^{iv}  }{e^{2n\pi  e^{2iv}}+1} )\diff v-\nsum{k=0}{n-1}\frac{\sqrt{2\pi}}{\sqrt{2k+1}})$$
$$=\lim_{n\to\infty}(\sqrt{n} \teg{0}{\pi } \ima(\frac{e^{iv}  \sqrt{2\pi}}{e^{2n\pi  e^{2iv}}+1} )\diff v-\nsum{k=0}{n-1}\frac{\sqrt{2\pi}}{\sqrt{2k+1}})$$
Although the suprema can not be directly applied to $\nsum{k=0}{n-1}\frac{1}{\sqrt{2k+1}}$ since it does not converge alone, it can be applied by moving the limit of $n$ to $\infty$ inside the integral such as :
\vskip-0.5cm
\begin{equation}\label{final yes hn real pn monotone}\boxed{\frac{1}{\sqrt{\pi}}\teg{-\infty}{\infty}\frac{\diff x}{e^{x^2}+1} =\lim_{n\to\infty}(\sqrt{2n} \teg{0}{\pi } \lim_{n\to\infty} \ima(\frac{e^{iv}  }{e^{2n\pi  e^{2iv}}+1} )\diff v-\sqrt{2}\nsum{k=0}{n-1}\frac{1}{\sqrt{2k+1}})}
\end{equation}
Because limits and imaginary part commute :
$$\teg{0}{\pi } \lim_{n\to\infty} \ima(\frac{e^{iv}  }{e^{2n\pi  e^{2iv}}+1} )\diff v=\teg{0}{\pi }  \ima(\lim_{n\to\infty}\frac{e^{iv}  }{e^{2n\pi  e^{2iv}}+1} )\diff v$$
$$=\teg{0}{\pi }  \ima(\lim_{n\to\infty}\frac{e^{iv}  }{e^{2n\pi  \cos(2v)}e^{2in\pi\sin(2v)}+1} )\diff v$$
For $0\le v<\pi/4$ and $3\pi/4 <v\le \pi$, we have that $0<\cos(2v)$. Because $|e^{2in\pi\sin(2v)}|=1$, no matter what value $\lim_{n\to\infty}e^{2in\pi\sin(2v)}$  goes to, the integrated function goes to zero on that interval because $\lim_{n\to\infty}e^{2n\pi  \cos(2v)}$ goes to infinity.
$$\teg{0}{\pi } \lim_{n\to\infty} \ima(\frac{e^{iv}  }{e^{2n\pi  e^{2iv}}+1} )\diff v\diff v=\teg{\pi/4}{3\pi/4 }  \ima(\lim_{n\to\infty}\frac{e^{iv}  }{e^{2n\pi  \cos(2v)}e^{2in\pi\sin(2v)}+1} )\diff v$$
For $\pi/4<v<3\pi/4$, we have that $\cos(2v)<0$, hence, $\lim_{n\to\infty}e^{2n\pi  \cos(2v)}$ goes to zero, which makes the function go to $e^{iv}$
\begin{equation}\label{eiv}
  \boxed{
    \teg{0}{\pi } \lim_{n\to\infty} \ima(\frac{e^{iv}  }{e^{2n\pi  e^{2iv}}+1} )\diff v=\teg{\pi/4}{3\pi/4 }  \ima(\frac{e^{iv}  }{0+1} )\diff v=\teg{\pi/4}{3\pi/4 }  \sin(v)\diff v=\sqrt{2}}
\end{equation}
\vskip-0.5cm
The singular points at $v=\pi/4, 3\pi/4$ can be ignored because it makes the limit go to $1/2$, which is a finite singular point.
\newline
Applying this to  \ref{final yes hn real pn monotone}:
\vskip-0.4cm
\begin{equation}\label{limitapplied yes v} \boxed{\frac{1}{\sqrt{\pi}}\teg{-\infty}{\infty}\frac{1}{e^{x^2}+1}\diff x =\lim_{n\to\infty}(2\sqrt{n}-\sqrt{2}\nsum{k=0}{n-1}\frac{1}{\sqrt{2k+1}})}
\end{equation}
We set :
\vskip-1cm
\begin{equation}\label{an first}
\boxed{A_n=2\sqrt{n}-\sqrt{2}\nsum{k=0}{n-1}\frac{1}{\sqrt{2k+1}}}
\end{equation}
\vskip-0.3cm
\begin{equation}\label{tn first}
\boxed{T_n=2\sqrt{n}-\nsum{k=1}{n}\frac{1}{\sqrt{k}}}
\end{equation}
We have :
\vskip-0.7cm
$$A_n-T_n=2\sqrt{n}-\sqrt{2}\nsum{k=0}{n-1}\frac{1}{\sqrt{2k+1}}-2\sqrt{n}+\nsum{k=1}{n}\frac{1}{\sqrt{k}}=\sqrt{2}(\nsum{k=1}{n}\frac{1}{\sqrt{2k}}-\nsum{k=0}{n-1}\frac{1}{\sqrt{2k+1}})$$
\begin{equation}\label{tn second}
  A_n-T_n=\sqrt{2}\nsum{k=1}{2n}\frac{(-1)^k}{\sqrt{k}}\longleftrightarrow \boxed{T_n=A_n+\sqrt{2}\nsum{k=1}{2n}\frac{(-1)^{k+1}}{\sqrt{k}} }
\end{equation}
We can also have, from \ref{an first} \ref{tn first} :
\vskip-0.5cm
$$A_n+T_n=4\sqrt{n}-\sqrt{2}\nsum{k=0}{n-1}\frac{1}{\sqrt{2k+1}}-\sqrt{2}\nsum{k=1}{n}\frac{1}{\sqrt{2k}}=4\sqrt{n}-\sqrt{2}(\nsum{k=0}{n-1}\frac{1}{\sqrt{2k+1}}+\nsum{k=1}{n}\frac{1}{\sqrt{2k}})$$
\begin{equation}\label{tn third}
\boxed{A_n+T_n=\sqrt{2}(2\sqrt{2n}-\nsum{k=1}{2n}\frac{1}{\sqrt{k}})=T_{2n}\sqrt{2}}
\end{equation}
Combining \ref{tn third} with \ref{tn second} :
\vskip-0.4cm
$$A_n+A_n+\sqrt{2}\nsum{k=1}{2n}\frac{(-1)^{k+1}}{\sqrt{k}}=(A_{2n}+\sqrt{2}\nsum{k=1}{4n}\frac{(-1)^{k+1}}{\sqrt{k}})\sqrt{2}$$
\begin{equation}\label{second an}
 \boxed{\sqrt{2}A_{2n}-2A_n=\sqrt{2}\nsum{k=1}{2n}\frac{(-1)^{k+1}}{\sqrt{k}}-2\nsum{k=1}{4n}\frac{(-1)^{k+1}}{\sqrt{k}}}
\end{equation}
We can have that, from \ref{an first}:
\vskip-1cm
$$A_{n+1}-A_n=2\sqrt{n+1}-2\sqrt{n}-\nsum{k=0}{n}\frac{1}{\sqrt{2k+1}}+\nsum{k=0}{n-1}\frac{1}{\sqrt{2k+1}}=2\sqrt{n+1}-2\sqrt{n}-\frac{1}{\sqrt{2n+1}}$$
We see that the function $2\sqrt{n+1}-2\sqrt{n}-\frac{1}{\sqrt{2n+1}}$ is always positive for all integers including zero $n$, implying $0<A_{n+1}-A_n$, $A_n<A_{n+1}$, hence $A_n$ strictly increases. We know that the limit $\lim_{n\to\infty} A_n=\frac{1}{\sqrt{\pi}}\teg{-\infty}{\infty}\frac{1}{e^{x^2}+1}\diff x$ must converge since the integral is between $0$ and $\frac{\pi}{\sqrt{2}}$ \ref{bounded}, implying, because $A_n$ is strictly increasing, $0\le A_n <\lim_{n\to\infty} A_n\le\sqrt{\pi/2}$. This bounds $A_n$, satisfying the monotone convergence theorem for $A_n$.
\newline
We then set $f_n=\nsum{k=1}{2n}\frac{(-1)^{k+1}}{\sqrt{k}}$
$$ f_{n+1}-f_n=\nsum{k=1}{2n+2}\frac{(-1)^{k+1}}{\sqrt{k}}-\nsum{k=1}{2n}\frac{(-1)^{k+1}}{\sqrt{k}}=\frac{1}{\sqrt{2n+1}}-\frac{1}{\sqrt{2n+2}}$$
We can also see that $\frac{1}{\sqrt{2n+1}}-\frac{1}{\sqrt{2n+2}}$ is always positive, implying $f_n$ is increasing. Furthermore, the convergence of $\lim_{n\to\infty} f_n$ is ensured from the alternating series test. We can also have that $\nsum{k=1}{n}\frac{(-1)^{k+1}}{\sqrt{k}}$ is at it's highest when $n=1$, bounding $f_n$ by 1, the monotone convergence theorem is then applicable to $f_n$ and $A_n$ \cite{moncontheo}. Hence, setting $\lim_{n\to\infty} A_n=A$, from \ref{second an} :
\vskip-1cm
$$\sqrt{2}A-2A=\lim_{n\to\infty}(\sqrt{2}A_{2n}-2A_n)=\lim_{n\to\infty}(\sqrt{2}\nsum{k=1}{2n}\frac{(-1)^{k+1}}{\sqrt{k}}-2\nsum{k=1}{4n}\frac{(-1)^{k+1}}{\sqrt{k}})$$
\vskip-0.5cm
$$(\sqrt{2}-2)A=(\sqrt{2}-2)\nsum{k=1}{\infty}\frac{(-1)^{k+1}}{\sqrt{k}} $$
$$A=\lim_{n\to\infty}( 2\sqrt{n}-\sqrt{2}\nsum{k=0}{n-1}\frac{1}{\sqrt{2k+1}})=\nsum{k=1}{\infty}\frac{(-1)^{k+1}}{\sqrt{k}}$$
From \ref{limitapplied yes v} :
\vskip-1cm
\begin{equation}\label{yessum}
    \boxed{\sqrt{\pi}\nsum{k=1}{\infty} \frac{(-1)^{k+1}}{\sqrt{k}}=\teg{-\infty}{\infty}\frac{1}{e^{x^2}+1} \diff x}
\end{equation}
We have that, from finite sums :
\vskip-0.4cm
$$\nsum{w=1}{n} (-1)^{w+1} x^w = \frac{x-(-1)^{n}x^{n+1} }{1+x}$$
Setting $x \to e^{-x^2}$
\vskip-0.4cm
$$\nsum{w=1}{n} (-1)^{w+1} e^{-wx^2} = \frac{e^{-x^2}-(-1)^{n} e^{-x^2(n+1)} }{1+e^{-x^2}} = \frac{1-(-1)^{n} e^{-x^2n} }{e^{x^2}+1}$$
$$ \frac{(-1)^{n} e^{-x^2n} }{e^{x^2}+1}+\nsum{w=1}{n} (-1)^{w+1} e^{-wx^2} = \frac{1}{e^{x^2}+1}$$
Applying this to \ref{yessum} :
\vskip-0.3cm
$$\sqrt{\pi}\nsum{k=1}{\infty} \frac{(-1)^{k+1}}{\sqrt{k}}=\teg{-\infty}{\infty} \frac{(-1)^{n} e^{-x^2n} }{e^{x^2}+1}+\nsum{w=1}{n} (-1)^{w+1} e^{-wx^2} \diff x$$
$$=\lim_{b\to\infty}(\teg{-b}{b} \frac{(-1)^{n} e^{-x^2n} }{e^{x^2}+1}\diff x+\teg{-b}{b} \nsum{w=1}{n} (-1)^{w+1} e^{-wx^2} \diff x)$$
Because finite sums interchange with integration :
\vskip-0.3cm
$$\sqrt{\pi}\nsum{k=1}{\infty} \frac{(-1)^{k+1}}{\sqrt{k}}=\lim_{b\to\infty}(\teg{-b}{b} \frac{(-1)^{n} e^{-x^2n} }{e^{x^2}+1}\diff x+\nsum{w=1}{n}\teg{-b}{b}  (-1)^{w+1} e^{-wx^2} \diff x)$$
\newline
Substituting $x\sqrt{w} =y$, $w x^2=y^2$, $\diff x=\diff y/\sqrt{w}$ on the right integral : 
\vskip-1cm
$$\sqrt{\pi}\nsum{k=1}{\infty} \frac{(-1)^{k+1}}{\sqrt{k}}=\lim_{b\to\infty}(\teg{-b}{b} \frac{(-1)^{n} e^{-x^2n} }{e^{x^2}+1}\diff x+\nsum{w=1}{n}\teg{-b \sqrt{w}}{b\sqrt{w}}  (-1)^{w+1} e^{-y^2} \frac{\diff y}{\sqrt{w}})$$
\vskip-0.3cm
$$=\lim_{b\to\infty}(\teg{-b}{b} \frac{(-1)^{n} e^{-x^2n} }{e^{x^2}+1}\diff x+\nsum{w=1}{n}\frac{(-1)^{w+1}}{\sqrt{w}} \teg{-b \sqrt{w}}{b\sqrt{w}}   e^{-y^2} \diff y)$$
Because $0\le e^{-x^2}\le 1/(x^2+1)$ bounds the convergence of $\teg{-\infty}{\infty} e^{-x^2}\diff x$ between $0$ and $\pi$, we see that $\lim_{b\to\infty} \teg{-b \sqrt{w}}{b\sqrt{w}}   e^{-y^2} \diff y$ converges between $0$ and $\pi$, hence we can separate the limit :
\vskip-0.7cm
$$\sqrt{\pi}\nsum{k=1}{\infty} \frac{(-1)^{k+1}}{\sqrt{k}}=\teg{-\infty}{\infty} \frac{(-1)^{n} e^{-x^2n} }{e^{x^2}+1}\diff x+\nsum{w=1}{n}\frac{(-1)^{w+1}}{\sqrt{w}} \teg{-\infty}{\infty}   e^{-y^2} \diff y$$
Setting $ B_n=(-1)^n \teg{-\infty}{\infty} \frac{ e^{-x^2n} }{e^{x^2}+1}\diff x$, because $\frac{1}{e^{x^2}+1}\le \frac{1}{x^2+2} \le \frac{1}{2}$ :
\vskip-0.2cm
$$ 0\le |B_n|\le \frac{1}{2}\teg{-\infty}{\infty}  \frac{e^{-x^2n} }{x^2+2}\diff x\le \frac{1}{2}\teg{-\infty}{\infty}  e^{-x^2n} \diff x=\frac{1}{2\sqrt{n} }\teg{-\infty}{\infty}e^{-x^2}\diff x$$
We already know $\teg{-\infty}{\infty}e^{-x^2}\diff x$ converges to a constant between $0$ and $\pi$, hence, letting $n$ go to $\infty$ :
\vskip-1cm
$$ 0\le\lim_{n\to\infty} |B_n|\le \lim_{n\to\infty}\frac{1}{2\sqrt{n} }\teg{-\infty}{\infty}e^{-x^2}\diff x=0$$
$0\le\lim_{n\to\infty} |B_n|\le 0$ gives us $0=\lim_{n\to\infty} |B_n|$, because limits and absolute value commute, $0=|\lim_{n\to\infty} B_n|$, which implies  $0=\lim_{n\to\infty} B_n$. If we let $n$ go to $\infty$ :
\vskip-0.4cm
$$\sqrt{\pi}\nsum{k=1}{\infty} \frac{(-1)^{k+1}}{\sqrt{k}}= \lim_{n\to\infty}(B_n+\nsum{w=1}{n}\frac{(-1)^{w+1}}{\sqrt{w}} \teg{-\infty}{\infty}  e^{-y^2} \diff y)=\nsum{w=1}{\infty}\frac{(-1)^{w+1}}{\sqrt{w}} \teg{-\infty}{\infty}  e^{-y^2} \diff y$$
$$\sqrt{\pi}\nsum{k=1}{\infty} \frac{(-1)^{k+1}}{\sqrt{k}}=\nsum{w=1}{\infty} \frac{(-1)^{w+1}}{\sqrt{w}}  \teg{-\infty}{\infty}  e^{-y^2} \diff y$$
$$ {\sqrt{\pi}=\teg{-\infty}{\infty}e^{-y^2}\diff y}$$
\begin{center}
  \textup{“} Cauchy's theorem cannot be employed to evaluate all definite integrals; thus $\int_{0}^{\infty} e^{-x^2}\diff x$ has not been evaluated except by other methods.\textup{”}, Watson, 1914, page 79, \cite{watson}
  \end{center}

\section{Additional formulas and link to Riemann Zeta}
There are a few extra formulas for solving the Gaussian integral that can be derived from this paper.
\newline
\newline
From \ref{norescont}, if we set $H_n$ such that it is lesser than the absolute value of the first pole, \ref{poles}, $H_n <\sqrt{\pi}$, we get that there would be no poles inside our contour hence no residue, thus :
$$  {\teg{-H}{H}\frac{1}{e^{x^2}+1} \diff x=-\teg{0}{\pi }\frac{i H  e^{iv} }{e^{H^2 e^{2iv}}+1} \diff v   \;\;\;\;\;\;\;\;\;\;\;\;( 0\le H<\sqrt{\pi})}$$
With this, we can get that the integral is actually always real as long as $H$ is real, hence:
\vskip-0.4cm
$$ 0=\rea(\teg{0}{\pi }\frac{ e^{iv} }{e^{H e^{2iv}}+1} \diff v) \;\;\;\;\;\;\;\;\;\;\;\;\;\;  (H \in \mathbb{R} )$$
With this and \ref{eiv}, we can get :
\vskip-1cm
$$ {i\sqrt{2}=\lim_{n\to\infty}\teg{0}{\pi} \frac{e^{iv} }{e^{2n\pi e^{2iv}}+1 }\diff v}$$
From \cite{finitezeta} :
\vskip-1cm
$$\zeta(x)=\lim_{n\to\infty} (\frac{n^{1-x}}{1-x}-\nsum{k=1}{n}\frac{1}{k^x})    (0<\rea(x))$$
$$\zeta(1/2)=\lim_{n\to\infty} (2\sqrt{n}-\nsum{k=1}{n}\frac{1}{\sqrt{k}})   $$
If we take \ref{limitapplied yes v} \ref{tn first} and \ref{tn third} and apply a limit, we would  get :
$$ {\sqrt{\pi}\zeta(1/2)(1-\sqrt{2})=\teg{-\infty}{\infty}\frac{1}{e^{x^2}+1}\diff x}$$
Furthermore, if we instead take \ref{yessum} and \cite{alternatingzeta}
$$\zeta(x)=\frac{1}{1-2^{1-x}}\nsum{k=1}{\infty}\frac{(-1)^{k+1}}{k^x}    (0<\rea(x))$$
$$\zeta(1/2)(1-\sqrt{2})=\nsum{k=1}{\infty}\frac{(-1)^{k+1}}{\sqrt{k}}  $$
We would also get :
\vskip-1cm
$$ \sqrt{\pi}\zeta(1/2)(1-\sqrt{2})=\teg{-\infty}{\infty}\frac{1}{e^{x^2}+1}\diff x$$
\newline
Which incites the proposition that this paper is actually a proof of the Riemann Zeta functional equation centered at $s=1/2$. If one were to prove it in a way similar to this paper, one would need to :
\newline
\newline
Define a function as :
\vskip -1cm
$$f(y)=\teg{0}{\infty}\frac{\diff x}{e^{x^y}+1}$$
Use finite sums like below and the same approach used before 
$$\nsum{w=1}{n} (-1)^{w+1} x^w = \frac{x-(-1)^{n}x^{n+1} }{1+x}$$
To get 
\vskip -1cm
$$\frac{1}{y}\nsum{w=1}{\infty}(-1)^{w+1}w^{-1/y} \teg{0}{\infty} e^{-k }k^{1/y-1}  \diff k = \nsum{w=1}{\infty}(-1)^{w+1} \teg{0}{\infty} e^{-w x^y }\diff x = f(y)$$
By \cite{alternatingzeta} and the definition of the gamma function :
$$\frac{(1-2^{1-1/y})}{y}\zeta(1/y)\Gamma(1/y) = f(y)$$
And then, using contour integration, obtaining a formula for $f(y)$ that is similar to \ref{limitapplied yes v}, doing the same manipulation tricks of \ref{an first} and \ref{tn first}.

\section{References}
\bibliography{refs} 
\bibliographystyle{plainurl}

\end{document}